\documentclass[conference]{IEEEtran}
 \usepackage{float}
\usepackage{epsfig}
\usepackage{mathptmx}
\usepackage{amssymb}
\usepackage{latexsym,graphicx}
\usepackage[usenames]{color}
\usepackage{amsmath,amssymb}
\usepackage{verbatim}
\newtheorem{theorem}{Theorem}
\newtheorem{definition}[theorem]{Definition}

\newtheorem{proposition}[theorem]{Proposition}
\newtheorem{remark}[theorem]{Remark}

\newtheorem{lemma}[theorem]{Lemma}

\newcommand{\erre}{\mbox{$\mathbb{R}$}}

\newcommand{\fine}{~\hspace*{\fill}{$\Box$}}

\begin{document}
\title{Order-type Henstock and McShane integrals in Banach lattice setting}
\author{\IEEEauthorblockN{ Domenico Candeloro\IEEEauthorrefmark{1} and 
Anna Rita Sambucini\IEEEauthorrefmark{7}}
\IEEEauthorblockA{\IEEEauthorrefmark{1}University of Perugia--Dipartimento di Matematica e Informatica, via Vanvitelli, 1 I-06123 Perugia (Italy) }
\IEEEauthorblockA{\IEEEauthorrefmark{7}University of Perugia--Dipartimento di Matematica e Informatica, via Vanvitelli, 1 I-06123 Perugia (Italy) }
\IEEEauthorblockA{E-mails:domenico.candeloro@unipg.it;anna.sambucini@unipg.it}}
\maketitle
\date{July $7^{\rm th}$, 2014}
\begin{abstract}  
Henstock-type integrals are studied for functions defined in a compact metric space $T$ endowed with a regular $\sigma$-additive measure $\mu$, and taking values in a Banach lattice $X$. In particular, the space $[0,1]$ with the usual Lebesgue measure is considered. \\
The norm- and the order-type integral are compared and interesting results are obtained when $X$ is an $L$-space.\\
\end{abstract}

\normalsize

\section{Introduction}
It is well-known that many different notions of integral were introduced in the last century, for real-valued functions, in order to generalize the Riemann one. An exhaustive discussion of the various definitions for real-valued functions can be found in \cite{KS}, where the Henstock-Kurzweil,
the Lebesgue, and
the McShane notions have been finely compared, and the equivalence between 
the McShane and
the Lebesgue integrals 
is
clearly described.   
The situation changes deeply in the case of Banach space-valued functions: in this case, it is well-known that the stronger type of integral is the Bochner one, which implies 
the Birkhoff integrability, which in turn is stronger than 
the McShane and the latter is stronger than both Henstock and Pettis integrals. The wide literature in this topic witnesses the great interest for these problems: see for example 
\cite{bbs,BS2004,BS2011,BS2004-nitra,bcs2011,bvg,cao,f1995,DPMILL,DPVM,dpp,dr,f1994a,f1994b,riecan,r2009a,r2009b,r2005,SS}. 
Alternative notions of integrals have also been given, for various applications (see e.g. \cite{ims1998a,ims1998b,is1998,LW2}).
 Subsequently, 
 the notions of order-type integrals have also been introduced and studied, for functions taking their values in ordered vector spaces, and in Banach lattices in particular: see
\cite{bvl,fvol3,mn,bd2014,bdp2012,bms,csmed,mimmoroma,dallas}.
 Also the multivalued case has been intensively studied, for both types of convergences: see for example \cite{cascales2007b,bcs2014}.

In this research mainly the differences between norm- and order-type McShane integrals 
have been 
 discussed, for single-valued functions taking values in a Banach lattice with order-continuous norm. 

However this note is just an anticipation of a forthcoming  extended paper, with further results and more detailed proofs.

After a section of preliminaries, in the third one  norm and order-type integrals are compared, showing 
a first striking difference: order integrals in general do not respect almost everywhere equality, except for order-bounded functions. Another interesting difference is that order integrals enjoy the so-called Henstock Lemma: this fact has interesting consequences in $L$-spaces,
where McShane order integrability is stronger than the Bochner (norm) one. 
In the fourth section  integrability in $[0,1]$ is discussed and it is proven that monotone mappings are McShane order-integrable, by using a similar procedure as in \cite{KS}.

\section{Preliminaries}
From now on, $T$ will denote a compact metric space, and $\mu:\mathcal{B}\to \erre^+_0$ any regular, nonatomic $\sigma$-additive measure on the $\sigma$-algebra $\mathcal{B}$ of Borel subsets of $T$. 

A \textit{gage} is any map $\gamma: T \rightarrow \mathbb{R}^+$.
A \textit{partition}  $\Pi$  of $[0,1]$ is a finite family $\Pi = \{ (E_i,t_i): i=1, \ldots,k \} $ of pairs such that the
 $E_i$ are pairwise disjoint sets whose union is $T$ and the points $t_i$ are called {\em tags}. If all tags satisfy the condition $t_i \in E_i$ then the partition is said to be of {\em Henstock} type, or a {\em Henstock partition}. Otherwise, it is said to be a {\em free} or {\em McShane} partition.

Given a gage $\gamma$,  a partition $\Pi$ is \textit{$\gamma$-fine} $(\Pi \prec \gamma)$ if
$d(w,t_i) < \gamma (t_i)$ for every $ w \in E_i$ and $i = 1, \ldots, k$.\\

Clearly, a gage $\gamma$ can  also be defined as a mapping associating with each point $t\in T$ an open ball centered at $t$: sometimes this concept will be used, without risk of confusion.\\

Let 
 $X$ be any Banach lattice with an order-continuous norm. 
For the sake of completeness  the main notions of
gauge integral are recalled here. \\

\begin{definition}\rm  \label{fnorm}
\rm A function $f:T\rightarrow X$ is \textit{norm-integrable}
 if there exists $J \in X$ such that, for every $\varepsilon > 0$
 there is a gage $\gamma : T \rightarrow \mathbb{R}^+$ such that for every $\gamma$-fine Henstock partition of $T$,

$\Pi=\{(E_i, t_i), i=1, \ldots, q \}$, it is:
\begin{eqnarray*}
\left\|\sigma(f, \Pi)- J \right\| \leq \varepsilon.
\end{eqnarray*}
(Here, as usual, the symbol $\sigma( f, \Pi)$ means
 $\sum_{i=1}^q  f(t_i) \mu(E_i)).$
In the case of integrability in the Henstock  sense, this will be denoted with \textit{{\rm H}-integrability} and the integral $J$ will be denoted with $H\int f d\mu$.\\
\end{definition}

\begin{remark}\label{mettipunt}\rm
It is worth noticing here that, in the previous definition, taking more generally {\em free} $\gamma$-fine partitions does not modify the concept introduced: in other words, the same integral is obtained if all free $\gamma$-fine partitions are allowed in the previous definition: see for example \cite[Proposition 2.3]{bcs2014}).
\end{remark}

Parallel to this definition, notions of {\em order-type} integral can be given, in accordance with the following
\begin{definition}\rm  \label{forder}
\rm A function $f:T\rightarrow X$ is \textit{order-integrable}
in the Henstock sense  if there exist $J \in X$,  an $(o)$-sequence $(b_n)_n$ in $X$ and a corresponding sequence $(\gamma_n)_n$ of gages, such that
for every $n$ and every $\gamma_n$-fine Henstock partition of $T$,
$\Pi=\{(E_i, t_i), i=1, \ldots, q \}$, one has
\begin{eqnarray*}
\left|\sigma(f, \Pi)- J \right| \leq b_n.
\end{eqnarray*}
In the case of integrability in the Henstock  sense, this will be denoted with \textit{{\rm (oH)}-integrability} and the integral $J$ will be denoted with $(oH)\int f$.\\
\end{definition}

It is obvious that also in this case there is no difference in taking all {\em free} $\gamma_n$-fine partitions in the previous definition.
Thanks to the order continuity of the norm in $X$, it is easy to see that any (oH)-integrable map $f$ is also H-integrable and the integrals coincide.
Furthermore, it can be  observed that (oH)-integrability of a function $f:T\to X$ implies also Pettis integrability, thanks to well-known results concerning the McShane norm-integral: see \cite[Theorem 8]{f1994a}.
However, later it will be shown that the (H)- and the (oH)-integrability are not equivalent, in general.

\section{Comparison between Norm and Order integral}
There are deep differences between order-type and norm-type integrals.
A first remarkable fact is that, in general, almost equal functions can behave in different ways with respect to the (oH)-integral, as it was proven in \cite[Example 2.8]{bms}, where a function $f:[0,1]\to c_{00}$ is 
given
with the following properties: $f$ is almost everywhere null (with respect to the Lebesgue measure) and $f$ is not (oH)-integrable. So, this function is almost everywhere equal to 0, hence is Bochner-integrable, but not  (oH)-integrable
or order-bounded.
Indeed, for order-bounded functions, the situation is better, as shown in the next Proposition.

\begin{proposition}
Let $f,g:T\to X$ be two bounded maps, such that $f=g$ $\mu$-almost everywhere.  Then, $f$ is {\rm (oH)}-integrable if and only if $g$ is, and the integral is the same.
\end{proposition}
{\bf Proof:}\ Let $M$ be any majorant for $|f|$ and $|g|$, and assume that $f$ is 
(oH)-integrable, with integral $J$. Let  $(b_n)_n$ and $(\gamma_n)_n$be  the sequences related to (oH)-integrability of $f$. In order to show integrability of $g$, fix $n$, and pick any open set $A_n\subset T$, with $\mu(A_n)< n^{-1}$ and 
$A_n \supset N:=\{t\in T:f(t)\neq g(t)\}$.
Now, for each element $u\in N$ let $\delta_n(u)$ be any open set containing $u$ and contained in $A_n$: then define $\gamma'_n(t)=\gamma_n(t)$ when $t\notin N$, while $\gamma_n'(u)=\gamma_n(u)\cap \delta_n(u)$ when $u\in N$.\\
Now, fix any tagged $\gamma'_n$-fine partition $\Pi$, ($\Pi:=(E_i,\tau_i)_i)$ and observe that, whenever
 the tag $\tau_i$ belongs to $N$, then $E_i\subset A_n$. 
So, it follows easily that
\begin{eqnarray*}
\sup \left\{
\sum_{\tau_i\in N} |f(\tau_i)|\mu(E_i),
\sum_{\tau_i\in N} |g(\tau_i)|\mu(E_i) \right\} \leq \frac{M}{n},
\end{eqnarray*}
while
$$\sum_{\tau_i\notin N}f(\tau_i)\mu(E_i)=\sum_{\tau_i\notin N}g(\tau_i)\mu(E_i),$$
and finally
\begin{eqnarray*}
&& |\sigma(g,\Pi)-J| \leq\\ 
&&\leq |\sigma(f,\Pi)-J|+\sum_{\tau_i\in N}|f(\tau_i)|\mu(E_i)+\\
&&+\sum_{\tau_i\in N}|g(\tau_i)|\mu(E_i)\leq b_n+2\frac{M}{n}.
\end{eqnarray*}
This clearly proves that $g$ is (oH)-integrable with integral $J$. 
Finally, interchanging the role of $f$ and $g$, the 
assertion follows.\ \fine
\\

\medskip

Another interesting difference is in the validity of the so-called Henstock Lemma: indeed, 
 (oH)-integrability yields this result, contrarily to the case of the norm-integral.\\

A  Cauchy-type criterion is stated first in order to prove the existence of the (oH)-integral. The proof is straightforward. 

\begin{theorem}\label{ordercauchy}
Let $f:T\rightarrow X$ be any mapping. Then $f$ is {\rm (oH)}-integrable  if and only if there exist an $(o)$-sequence $(b_n)_n$ and a corresponding sequence $(\gamma_n)_n$ of gages, such that for every $n$, as soon as $\Pi, \Pi'$ are two $\gamma_n$-fine Henstock  partitions, the following holds:
\begin{eqnarray}\label{cauchyprimo}
|\sigma(f,\Pi)-\sigma(f,\Pi')|\leq b_n
\end{eqnarray}
\end{theorem}
Now, a Henstock-type lemma is stated, for the (oH)-integral.
The proof is similar to that of  \cite[Theorem 1.4]{mimmoroma}, and follows from the Cauchy criterion. However, 
some details are also given here.\\

\begin{proposition}\label{henlemma}
Let $f:T\to X$ be any {\rm (oH)}-integrable  function. Then, there exist an $(o)$-sequence $(b_n)_n$ and a corresponding sequence $(\gamma_n)_n$ of gages, such that, for every $n$ and every $\gamma_n$-fine Henstock  partition $\Pi$ it 
is
$$\sum_{E\in \Pi}Ob_n(f,E)\leq b_n,$$
where $$Ob_n(E)=\sup_{\Pi'_E,\Pi''_E}\{|\hskip-2mm\sum_{F''\in \Pi''_E}f(\tau_{F''})\mu(F'') - \hskip-2mm\sum_{F'\in \Pi'_E}f(\tau_{F'})\mu(F')|\},$$
and $\Pi'_E, \Pi''_E$ run along all $\gamma_n$-fine Henstock partitions of $E$.
\end{proposition}
{\bf Proof:}\  First observe that, thanks to the Cauchy criterion, an $(o)$-sequence $(b_n)_n$ exists, together with a corresponding sequence of gages $(\gamma_n)_n$, such that
\begin{eqnarray}\label{primocauchy}
|\sum_{F'\in \Pi'}f(\tau_{F'})\mu(F')-\sum_{F''\in \Pi''}f(\tau_{F''})\mu(F'')|\leq b_n
\end{eqnarray}
(with obvious meaning of symbols)
holds, for all $\gamma_n$-fine partitions $\Pi', \ \Pi''$. Now, take any $\gamma_n$-fine partition $\Pi$, and, for each element $E$ of $\Pi$, consider two arbitrary subpartitions $\Pi'_E$ and $\Pi''_E$.
Then, taking the {\em union} of the subpartitions $\Pi'_E$ as $E$ varies, and making the same operation with the subpartitions $\Pi''_E$,  two $\gamma_n$-fine partitions of $T$ are obtained, for which (\ref{primocauchy}) holds true. From  (\ref{primocauchy}),  obviously it follows
\begin{eqnarray}\label{secondocauchy}
\sum_{F'\in \Pi'}f(\tau_{F'})\mu(F')- \sum_{F''\in \Pi''}f(\tau_{F''})\mu(F'')\leq b_n.
\end{eqnarray}
Now,   let $E_1$ be the first element of $\Pi$.
In the summation at left-hand side, fix all the $F's$ and the $F''s$ that are not contained in $E_1$. 
Taking the supremum when the remaining $F's$ and $F''s$  vary in all possible ways, it follows
\begin{eqnarray*}
 Ob_n(f,E_1)+ \hskip-2mm \sum_{\substack{F'\in \Pi',\\ F'\not\subset E_1}}f(\tau_{F'})\mu(F') - \hskip-2mm
 \sum_{\substack{F''\in \Pi'',\\F''\not\subset E_1}}f(\tau_{F''})\mu(F'')\leq b_n.
\end{eqnarray*}
In the same fashion, fixed all the $F'$ and $F''$ that are not contained in the second subset of $\Pi$, (say $E_2$), and making the same operation, it follows 
\begin{eqnarray*}
\hskip-1cm&& Ob_n(f,E_1)+Ob_n(f,E_2)+\hskip-2mm \sum_{\substack{F'\in \Pi',\\ F'\not\subset E_1\cup E_2}}f(\tau_{F'})\mu(F') + \\
&&- \sum_{\substack{F''\in \Pi'',\\F''\not\subset E_1\cup E_2}}f(\tau_{F''})\mu(F'')\leq b_n
\end{eqnarray*}
Now, it is clear how to deduce the assertion. \fine\\

\begin{remark}\label{subinterval}\rm
A first consequence of the previous Proposition \ref{henlemma} is that any (oH)-integrable function $f$ is also integrable
 in the same sense in every measurable subset $A$. Indeed,
taking the same $(o)$-sequence $(b_n)_n$ and the same corresponding sequence $(\gamma_n)_n$ as for integrability of $f$, 
for each $n$ any  $\gamma_n$-fine partition of $A$ can be extended to a $\gamma_n$-fine partition of $T$ thanks to the Cousin Lemma, and so, for any two $\gamma_n$-fine partitions $\Pi, \ \Pi'$ of $A$, it follows
$$|\sigma(f,\Pi)-\sigma(f,\Pi')|\leq Ob_n(f,A)\leq b_n.$$
Then, the Cauchy criterion yields the conclusion.\\
\end{remark} 
\begin{remark}\label{additivity}\rm
By means of usual techniques, one also proves additivity of the integral:
namely
whenever  $f$ is integrable in $T$, and $A,B$ are two disjoint measurable subsets of $T$, then $\int f 1_{A\cup B}d\mu=\int_A f d\mu +\int_B f d\mu$.\\
\end{remark}

The following Theorem collects some easy consequences of  Proposition \ref{henlemma}.

\begin{theorem}\label{henstoc2}
Let $f:T\to X$ be any {\rm (oH)}-integrable function. Then there exist an $(o)$-sequence $(b_n)_n$ and a corresponding sequence $(\gamma_n)_n$ of gages, such that:
\begin{description}
\item[\ref{henstoc2}.1)]
 for every $n$ and every $\gamma_n$-fine partition $\Pi$ 
one has
$$\sum_{E\in \Pi}|f(\tau_E)\mu(E)-{\rm (oH)}\int_Ef d\mu|\leq b_n.$$
\item[\ref{henstoc2}.2)]
for every $n$ and every $\gamma_n$-fine partition $\Pi$ it holds
$$\sum_{E\in \Pi}|f(\tau_E)\mu(E)-f(\tau_E')\mu(E)|\leq b_n,$$
as soon as all the tags  satisfy the condition $E\subset \gamma_n(\tau_E')$ and $E\subset \gamma_n(\tau_E)$ for all $E$.
\end{description}
\end{theorem}
\begin{remark}\label{paracool}\rm
For further reference, 
observe
that in the theorem above all partitions may also be free, 
since, as  already noticed, the restriction $\tau_E\in E$ does not affect the results.
\end{remark}

A consequence of this theorem is that 
the (oH)-integrability of $f$ implies
the
 (oH)-integrability of $|f|$. 
\begin{theorem}\label{modulointegrabil}
If $f:T\to X$ is {\rm (oH)}-integrable, then also $|f|$ is.
\end{theorem}
{\bf Proof:} 
The Cauchy criterion is used in the following formulation:
there exist an $(o)$-sequence $(b_n)_n$ and a corresponding sequence $(\gamma_n)_n$ of gages, such that, for each $n$, as soon as $\Pi,\Pi'$ are $\gamma_n$-fine free partitions and $\Pi'$ is finer than $\Pi$, 
\begin{eqnarray}\label{cauchyy}
\big|\sum_{E\in \Pi}|f(\tau_E)|\mu(E)-\sum_{E'\in \Pi'}|f(\tau_{E'})|\mu(E')\big|\leq b_n.
\end{eqnarray}
Now, if $(b_n)_n$ and $(\gamma_n)_n$ are as in Theorem \ref{henstoc2}, and  $\Pi$ and $\Pi'$ are as  above, 
it is
\begin{eqnarray*}
&& \sum_{E\in \Pi}|f(\tau_E)|\mu(E)-\sum_{E'\in \Pi'}|f(\tau_{E'})|\mu(E') =\\
=&&\sum_{E\in \Pi}\sum_{\substack{E'\in \Pi',\\ E'\subset E}}(|f(\tau_E)|-|f(\tau_{E'})|)\mu(E')= \\
=&& \sum_{E'\in \Pi'}(|f(\tau^*_{E'})|-|f(\tau_{E'})|)\mu(E'),
\end{eqnarray*}
where the tags $\tau^*_{E'}$ coincide with $\tau_E$ whenever $E'\subset E$, $E\in \Pi$.
Therefore, a simple application of \ref{henstoc2}.2) (and Remark \ref{paracool}) leads to (\ref{cauchyy}) and the proof is finished. \fine 
\\

The last theorem can be compared with a previous result by Drewnowski and Wnuk, (\cite[Theorem 1]{DW}), where the Bochner integral is considered, for functions taking values in a Banach lattice $X$, and it is proven that the {\em modulus} of
 the indefinite Bochner integral of $f$ is precisely the indefinite integral of $|f|$. 

Here  a similar result for {\rm (oH)}-integrable mappings is stated, after introducing a new definition.
\begin{definition}\label{modulus}\rm
Let $f:T\to X$ be any {\rm (oH)}-integrable mapping, and set
$$\mu_f(A)={\rm (oH)}\int_A f d\mu$$
for all Borel sets $A\in \mathcal{B}$.
 Then $\mu_f$ is said to be the {\em indefinite integral} of $f$. The {\em modulus} of $\mu_f$, denoted by $|\mu_f|$, is defined for each $A\in \mathcal{B}$ as follows:
$$|\mu_f|(A)=\sup\{\sum_{B\in \pi}|\mu_f(B)|: \pi\in \Pi(A)\}$$
where $\Pi(A)$ is the family of all finite partitions of $A$.
(The {\em boundedness} of this quantity will be proven soon).\\
\end{definition}

Now 
the following theorem can be stated.\\
\begin{theorem}\label{parallelo}
Assume that $f:T\to X$ is {\rm (oH)}-integrable. Then 
one has
$$|\mu_f|= \mu_{|f|}.$$
\end{theorem}
{\bf Proof:} First of all, observe that $|f|$ is integrable too, thanks to Theorem \ref{modulointegrabil}. Since clearly
$$\left|{\rm(oH)}\int_B fd\mu \right|\leq {\rm (oH)}\int_B |f| d\mu$$
holds for every set $B\in \mathcal{A}$, it follows 
that
$$|\mu_f|\leq \mu_{|f|}.$$
This also shows that the modulus $|\mu_f|$ is bounded.
 Now, since $|\mu_f|$ and $\mu_{|f|}$ are additive, in order to obtain the reverse inequality, it will be sufficient to prove that $\mu_{|f|}(T)=|\mu_f|(T)$.
To this aim, 
the Henstock Lemma, \ref{henstoc2}, will be used, and in particular  \ref{henstoc2}.1). Let $(b_n)_n$ and $(\gamma_n)_n$  be an $(o)$-sequence and its corresponding sequence of gages, related to integrability of both $f$ and $|f|$.
 So, for every $n$ and every $\gamma_n$-fine partition $\pi\equiv (E_i,t_i)_i$ it holds
\begin{eqnarray*}
\hskip-.5cm
\mu_{|f|}(T)&-\hskip-.2cm&|\mu_f|(T) \leq \sum_i\left(\mu_{|f|}(E_i)-\left|{\rm (oH)}\int_{E_i}fd\mu \right|\right) \leq 
\\ \leq &&
 \sum_i\left (\mu_{|f|}(E_i)-|f(t_i)|\mu(E_i)\right)+\\
+&& 
\sum_i\left (|f(t_i)|\mu(E_i)-\left|{\rm (oH)}\int_{E_i}fd\mu\right|\right)\leq 2 b_n
\end{eqnarray*}
Since $(b_n)_n$ is an $(o)$-sequence, then $\mu_{|f|}(T)\leq|\mu_f|(T)$, and so clearly also $\mu_{|f|}(T)=|\mu_f|(T)$. This concludes the proof. \fine
\\

Another interesting consequence is concerned with $L$-spaces.  Recall that a Banach lattice $X$ is an $L$-space if its norm $\|\cdot\|$ satisfies
$$\|x+y\|=\|x\|+\|y\|$$
for all positive elements $x,y$ in $X$ (see also \cite{fvol3}).\\

The following definition, related with norm-integrability, is needed.\\
\begin{definition} \rm \cite[Definition 3]{DPVM}
 $f:T \to X$ is {\em variationally {\rm H} integrable} (in short {\rm vH}-integrable) if 
for every $\varepsilon>0$ there exists a gage $\gamma$ such that, for every $\gamma$-fine partition $\Pi\equiv(E,t_E)_E$ the following holds:
\begin{eqnarray}\label{variazio}
\sum_{E\in \Pi}\|f(t_E)\mu(E)-{\rm (H)}\int_E f d\mu\|\leq \varepsilon.
\end{eqnarray}
\end{definition}
For results on this setting see also \cite{DPVM}. \\
\begin{theorem}\label{Lspazio}
Let $f:T\to X$ be {\rm (oH)}-integrable, and assume that $X$ is an $L$-space. Then  $f$ is Bochner integrable.
\end{theorem}
{\bf Proof:}\ In order to prove  the Bochner integrability, it will suffice to show that 
$f$ is vH-integrable.
Indeed, by \cite[Theorem 2]{DPMILL},
the 
 variational integrability implies  
the 
Bochner integrability. 
The  H-integrability of $\|f\|$ will be proved first. 
In accordance with the previous Theorem \ref{henstoc2}, and with the same meanings of symbols, there exist an $(o)$-sequence $(b_n)_n$ and a corresponding sequence $(\gamma_n)_n$ of gages, such that, for every $n$ and every $\gamma_n$-fine partition $\Pi$ 
one has
$$\sum_{E\in \Pi}|f(\tau_E)\mu(E)-f(\tau_E')\mu(E)|\leq b_n.$$ 
Since the norm of $X$ is compatible with the order, 
then
$\lim_n\|b_n\|=0$. So, fix $\varepsilon>0$ and pick any integer $N$ such that $\|b_N\|\leq \varepsilon$.
Then, if $\Pi$ is any $\gamma_N$-fine partition, we have
$$\left\| \, \sum_{E\in \Pi}|f(\tau_E)\mu(E)-f(\tau_E')\mu(E)| \, \right\| \leq \|b_N \| \leq \varepsilon,$$
in accordance with  \ref{henstoc2}.2). Thanks to the particular nature of the norm $\|\cdot \|$, we deduce
that
\begin{eqnarray*}
\sum_{E\in \Pi}\left\| \, |f(\tau_E)\mu(E)-f(\tau_E')\mu(E)| \, \right\|\leq \|b_N\|\leq \varepsilon,
\end{eqnarray*}
and so
\begin{eqnarray}\label{normaepsilo}
\sum_{E\in \Pi} \, \big|\|f(\tau_E)\|-\|f(\tau_E')\|\big|\mu(E) \, \leq \|b_N\|\leq \varepsilon,\end{eqnarray}
as soon as $\Pi$ is $\gamma_N$-fine, both for the tags $\tau_E$ and for the tags $\tau_E'$. 
\\
Now, proceeding as in the proof of Theorem \ref{modulointegrabil}, it is not difficult to prove that $\|f\|$ satisfies the Cauchy criterion  for the Henstock integrability, and therefore it is integrable.
From (\ref{normaepsilo}), also thanks to \ref{henstoc2}.1), it is easy to deduce also (\ref{variazio}).
In conclusion, $f$ is variationally integrable, its norm is integrable, and then, from Pettis integrability, it follows also the Bochner integrability.
 \fine\\

\begin{remark}\label{HnonoH}\rm
The previous result can be used to show that 
the 
 (H)-integrability in general does not imply
the 
 (oH)-integrability: indeed, if $X$ is any infinite-dimensional Banach space, there exists a McShane (norm)-integrable map $f:[0,1]\to X$ that is not Bochner integrable (see \cite{SS}). In particular, when $X$ is an $L$-space (of infinite dimension), such function $f$ cannot be (oH)-integrable, in view of Theorem \ref{Lspazio}.
\end{remark}

\section{Integrability in $[0,1]$}
In this section 
the (oH)-integrability of functions $f:[0,1]\to X$ is studied, where
 $[0,1]$ is endowed with the usual Lebesgue measure $\lambda$ . \\

From now on, 
only {\em free} partitions consisting of subintervals of $[0,1]$ are considered, rather than arbitrary measurable subsets. Indeed, in \cite{f1995} it is proven that there is equivalence between the two types: though the proof there is related only to
 norm integrals, the technique is the same. For this reason, from now on 
the symbol
(oM)-integral will be used  rather than (oH)-integral. \\

A useful result, parallel to \cite[ Lemma 5.35]{KS}, is the following: \\
\begin{lemma}\label{come5.35}
Let $f:[0,1]\to X$ be any fixed function, and suppose that there exists an $(o)$-sequence $(b_n)_n$ such that, for every $n$ 
two {\rm (oM)}-integrable functions $g_1$ and $g_2$ can be found, with the same regulating $(o)$-sequence 
$(\beta_n)_n$, such that $g_1\leq f\leq g_2$ and ${\rm (oM)}\int g_2 d\lambda \leq {\rm (oM)}\int g_1 d\lambda+b_n.$
Then $f$ is {\rm (oM)}-integrable.
\end{lemma}
\medskip

The fact that  increasing functions are (oM)-integrable can be deduced similarly as in \cite[example 5.36]{KS}.
\begin{theorem}\label{mcmonotone}
Let $f:[0,1]\to X$ be increasing. Then $f$ is {\rm (oM)}-integrable. 
\end{theorem}
\medskip

\begin{remark}
In a similar way, one can prove 
the 
(oM)-integrability more generally for (order bounded) mappings that are Riemann-integrable in the order sense.
\end{remark}


\section{Conclusion}
In this paper   the notions of Henstock and McShane integrability for functions defined in a metric compact regular space and taking values in a Banach lattice with an order-continuous norm are investigated. Both the norm-type and the order-type integrals have been studied and compared. Though in general the order-type integral is stronger than the norm-one, 
in $M$-spaces the two notions coincide, while in $L$-spaces the order-type Henstock integral is indeed a Bochner one. Finally, the particular case of  functions defined in a real interval is considered, and it is proven that monotone mappings are always order-McShane integrable.

\vspace{1cm}
\noindent \textbf{Acknowledgment} The authors have been supported by
University of Perugia -- Department of Mathematics and Computer Sciences - Grant Nr 2010.011.0403, 
Prin "Metodi logici per il trattamento dell'informazione",  
Prin "Descartes"
 and by the Grant prot. U2014/000237 of GNAMPA - INDAM (Italy).


\end{document}